\documentclass[a4paper,11pt]{article}
\usepackage{amssymb,amsmath,amsthm}
\usepackage{a4wide}
\usepackage{pgf,tikz}
\usetikzlibrary{arrows}
\usepackage{epsfig}
\usepackage{color}

\usepackage{hyperref}

\usepackage{authblk}

\usepackage{subcaption}

\newcommand{\beq}{\begin{equation}}  
\newcommand{\eeq}{\end{equation}}  
\newcommand{\bea}{\begin{eqnarray}} 
\newcommand{\eea}{\end{eqnarray}}   
\newcommand{\bear}{\begin{array}}  
\newcommand{\eear}{\end{array}}

\newtheorem{thm}{Theorem}[section] 

\newtheorem{propn}[thm]{Proposition} 

\newtheorem{conje}[thm]{Conjecture}

\newenvironment{prf}{\trivlist \item [\hskip 
\labelsep {\bf Proof:}]\ignorespaces}{\qed \endtrivlist}

\theoremstyle{definition}
\newtheorem{definition}[thm]{Definition}
\newtheorem{remark}[thm]{Remark}
 
\newcommand{\T}{{\mathbb T}}

\newcommand{\Q}{{\mathbb Q}}
\newcommand{\Z}{{\mathbb Z}}
\newcommand{\C}{{\mathbb C}}
\newcommand{\R}{{\mathbb R}}
\newcommand{\Pro}{{\mathbb P}}

\newcommand{\rd}{\mathrm{d}}
\newcommand{\ri}{\mathrm{i}}

\newcommand\la{{\lambda}}

\newcommand\al{{\alpha}}

\newcommand\om{{\omega}}

\newcommand\si{{\sigma}}

\newcommand\bv{{\bf v}}
\newcommand\bx{{\bf x}}
\newcommand\bd{{\bf d}}

\newcommand\bmu{\boldsymbol\mu}
\newcommand\bo{\boldsymbol{1}}

\newcommand{\cA}{{\cal A}}
\newcommand{\cE}{{\cal E}}
\newcommand{\cO}{{\cal O}}
\newcommand{\cI}{{\cal I}}

\def\lc{\left\lfloor}   
\def\rc{\right\rfloor}

\begin{document}


\title{Growth of Mahler measure and algebraic entropy of dynamics with the Laurent property}
\author[1]{Andrew N. W. Hone} 
\affil[1]{School of Mathematics, Statistics \&  Actuarial Science, 
University of Kent, 
Canterbury CT2 7FS, U.K. 
}

\maketitle

\begin{abstract} 
We consider the growth rate of the Mahler measure in discrete dynamical systems with the Laurent property, 
and in cluster algebras, and compare this with other measures of growth. In particular, we formulate the conjecture 
that the growth rate of the logarithmic Mahler measure coincides with the algebraic entropy, which is defined 
in terms of degree growth. Evidence for this conjecture is provided by exact and numerical 
calculations of the Mahler measure for a family of Laurent polynomials generated by rank 2 cluster algebras, 
for a recurrence of third order related to the Markoff numbers, and for the Somos-4 recurrence. Also, for  
the sequence of Laurent polynomials associated with the Kronecker quiver (the cluster algebra of affine type $\tilde{A}_1)$
we prove a precise formula for the leading order asymptotics of the logarithmic Mahler measure, which grows linearly.
\end{abstract}

\section{Introduction, definitions and main conjecture} 

\setcounter{equation}{0}

 Given $P(x_1,\ldots,x_N)\in \C[x_1^{\pm 1},\ldots,x_N^{\pm 1}]\setminus\{0\}$, a non-zero Laurent polynomial in $N$ variables with complex coefficients,  
its logarithmic Mahler measure  is 
defined to be 
\beq\label{mm} 
m(P):=\int_{\T^N} \log | P(x_1,\ldots,x_n) |\, \rd \Omega, 
\eeq  
obtained by integrating over all points on the $N$-torus  with respect to the Haar measure 
\beq\label{haar} 
\rd \Omega =\frac{ \rd t_1\wedge\cdots \wedge \rd t_N}{(2\pi)^N} 
\eeq 
for $(x_1,\ldots,x_N)=(e^{\ri t_1},\ldots, e^{\ri t_N})\in\T^N$ with local coordinates $t_j\in (-\pi,\pi]$, 
which can be considered as the restriction to the real torus of the log-canonical meromorphic $N$-form 
\beq \label{vol} 
\rd \Omega = \frac{\rd x_1\wedge \cdots \wedge \rd x_N}{(2\pi\ri)^N\prod_{j=1}^N x_j} 
\eeq 
on $(\C^*)^N$ (or $\C^N$), with $\ri=\sqrt{-1}$.  
The name Mahler measure is also given to the quantity $M(P)=\exp m(P)$, but henceforth we will use the name without qualification to refer to the 
logarithmic version. For a polynomial $P(x)\in\C[x]$ in one variable of degree $d$, with leading coefficient $c\neq 0$ and roots 
$\al_j$, $j=1,\ldots, d$, the 
Mahler measure is found  explicitly by Jensen's formula to be 
\beq\label{jensen} 
m(P)=\log |c| +\sum_{j=1}^d [\log |\al_j|]_+ 
\eeq 
(here $[a]_+ =\max (a,0)$), 
and while the existence of the $N$-variable integral in (\ref{mm})  is  not immediately obvious, it is guaranteed by a result of Lawton \cite{lawton}, which shows 
that this reduces to the case of one variable, in the sense that 
\beq\label{limk}
m\big(P(x_1,\ldots,x_N)\big)=\lim_{k_2,\ldots, k_N\to \infty} m\big( P(x,x^{k_2},\ldots, x^{k_N})\big). 
\eeq    
In particular, this implies that $m(P)\geq 0$ for Laurent polynomials with integer coefficients. There are many interesting arithmetical questions concerning 
Mahler measure, like Lehmer's problem about the smallest non-zero value of $m(P)$ for $P(x)\in\Z[x]$, as well as the fact 
that explicit evaluations of $m(P)$ for specific  multivariable polynomials $P$  can be expressed in terms of $L$-functions of zeta functions, or volumes of hyperbolic polyhedra \cite{boyd, ew, lalin}. The first such example 
in two variables is due to Smyth \cite{smyth}, who obtained the formula
\beq\label{smyth}
m(x_1+x_2+1)=L'(\chi_{-3},-1)
\eeq 
in terms of the $L$-function for the Dirichlet character $\chi_{-3}$.

Cluster algebras, introduced in \cite{fz1}, are a class of commutative algebras which, rather than being 
specified a priori, have distinguished sets of generators called clusters, which are produced recursively from an initial $N$-tuple 
of generators $\bx = (x_1,\ldots,x_N)$ by a process called mutation. 
One of the basic features of cluster algebras is the Laurent 
property: the $N$ elements in each cluster are Laurent polynomials in the initial generators, with integer coefficients. In fact, Zelevinsky once gave the following informal definition \cite{zel}: ``\textit{A cluster algebra is a machine for generating 
non-trivial Laurent polynomials}.''  More precisely, a coefficient-free cluster algebra 
$\cA=\cA(\bx, B)$  of rank $N$ 
is generated by starting from a seed $(\bx, B)$ consisting of an initial cluster $\bx$ of $N$ variables and 
a matrix $B\in  \mathrm{Mat}_N(\mathbb{Z})$ that is skew-symmetrizable, in the sense that there exists a diagonal 
matrix of positive integers $D$ such that $DB$ is skew-symmetric. 
Then  for each integer
$k\in [1,N]$ there is a  mutation ${\mu}_k$, an involution  which produces a new seed
$( {\bf x}',B')={\mu}_k ( {\bf x},B)$, where  $B'=(b_{ij}')$ with
\beq\label{matmut}
b_{ij}' =  \begin{cases}
 -b_{ij} &\text{if}  \,\,i=k \,\, \text{or} \,\, j=k , \\
 b_{ij}+\text{sgn}(b_{ik})[b_{ik}b_{kj}]_+ & \text{otherwise},
\end{cases}
\eeq
and ${\bf x}'=(x_j')$ with 
\beq \label{clustmut}
x_j'=
 \begin{cases}
x_k^{-1}\, \big( \prod_{j=1}^N x_j^{[b_{kj}]_+}+  \prod_{j=1}^N x_j^{[-b_{kj}]_+}  \big)  &\text{for}  \,\,j=k  \\
 \ x_j   &\text{for}  \,\, j \neq k.
\end{cases}
\eeq 
Now for an arbitrary  composition of mutations 
$ 
\bmu = \mu_{k_n}\circ \cdots \circ \mu_{k_1}
$ 
of any length $|\bmu |=n\geq 0$, 
the Laurent property says that every element of the cluster $\bmu (\bx)$ belongs to 
$\Z[x_1^{\pm 1},\ldots,x_N^{\pm 1}]$, and the  cluster algebra $\cA$ is defined to be the $\Z$-algebra generated by all 
cluster variables in all clusters, that is $\cA = \Z \big[\underset  { \bmu}{\cup}{\bmu}(\bx)\big]$. The cluster algebra $\cA$ 
is determined by the seed $(\bx,B)$ up to mutation equivalence, i.e.\ $\cA(\bx,B)$ is isomorphic to $\cA(\bx',B')$ for 
$(\bx,B)\sim (\bx',B')$, where the equivalence relation means that $(\bx',B')=\bmu(\bx,B)$ for some composition of mutations $\bmu$.

For discrete dynamical systems defined by iteration of a rational map $\varphi$ in dimension $N$, there are  various notions of entropy that are used to 
measure the growth of complexity of the iterates. Letting  $\deg \varphi^n$ denote the maximum of the degrees of the rational functions defining the 
components of $ \varphi^n(\bx)$, the algebraic entropy was defined in \cite{bv} as follows. 

\begin{definition}\label{algen}
The algebraic entropy of a rational map $\varphi$ is given by 
\beq\label{E}
\cE :=\underset {n\to\infty}{\lim} \frac{\log \deg\varphi^n}{n}. 
\eeq 
\end{definition}

Note that the existence of the above limit is guaranteed by the subadditive property of the log degrees and Fekete's lemma \cite{fekete}. This 
measure of growth was first used extensively to study certain birational transformations arising in statistical mechanics (see \cite{aa, bk} and references therein).
For birational maps $\Pro^2\to\Pro^2$ it is possible to perform blowups at singular points and reduce the computation of $\cE$ to a 
calculation in terms of the induced linear action of the pullback $\varphi^*$ on the cohomology of the resulting enlarged space \cite{df, tak}, but in dimension $N>2$ 
the problem is more difficult, because where the map has singularities it is not always possible to remove exceptional hypersurfaces from the dynamical system by blowups. In general, numerical computation of  
the algebraic entropy of $\varphi$ is extremely computationally intensive, because it requires exact calculation of the iterates of $\varphi$, and generically  $\deg\varphi^n$ grows exponentially with $n$. 

Due to the fact that matrix mutations (\ref{matmut}) change the exponents that appear in (\ref{clustmut}), and because there are $N$ possible mutations that can be applied to any cluster, generically a seed 
in a cluster algebra $\cA$ defines a set of evolutions on an $N$-regular tree rather than a single dynamical system (apart from the so-called bipartite belt when $N=2$) \cite{fz1}; 
nevertheless, one can define a corresponding notion of algebraic entropy, analogous to (\ref{E}),  via 
\beq\label{EA}
\cE(\cA):= \underset {n\to\infty}{\lim}\underset {|\bmu |=n}{\mathrm{sup}}\, \frac{\log \deg\bmu(\bx)}{n}, 
\eeq  
where $\deg\bmu(\bx)$ is the largest total degree of the elements in the cluster $\bmu(\bx)$, considered as rational functions of the variables in the initial cluster 
$\bx =(x_1,\ldots,x_N)$. This definition only depends on the mutation equivalence class of the seed $(\bx,B)$ defining $\cA$, as it should, because the algebraic entropy is invariant 
under birational transformations \cite{bv}.

Inspired by Vojta's dictionary between Nevanlinna theory and Diophantine approximation, Halburd proposed that when a rational map $\varphi$ is defined over $\Q$ (or a number field), instead of degree growth it is  easier (at least experimentally) to measure 
the growth of the logarithmic height $h\big(\varphi^n(\bx)\big)$ along an orbit $\cO$ with specific numerical values for the initial data $\bx$, where, for $\bx\in\Q^N$, 
$h(\bx)=\log H(\bx)$, with the height 
$H(\bx)$ being the maximum modulus of the numerators and denominators of the components of $\bx$ (written as fractions in lowest terms, with the convention that the number 0 has height 1). The fact that growth of (naive) height is a useful measure of complexity was already noted for birational transformations appearing in the context of statistical mechanics \cite{aa}, 
and the additional observations made by Halburd suggest that 
the following definition is worthwhile. 

\begin{definition}\label{ED} 
For a  rational map (\ref{phimap}), the Diophantine entropy of 
a particular non-singular orbit $\cO$ is  defined by  
\beq\label{dioph} 
\cE_H(\cO) :=\underset {n\to\infty}{\lim} \frac{\log h\big(\varphi^n(\bx)\big)}{n}, 
\eeq 
while the  Diophantine entropy of the map is taken to be 
$$\cE_H(\varphi):=
\mathrm{sup}\,\cE_H(\cO), $$ 
where the supremum is taken over all non-singular orbits $\cO$. 
\end{definition}

In the latter definition, 
where it is possible the supremum can be taken over all orbits of an extension 
of $\varphi$ to a regular morphism on an enlarged phase space (with the expectation being 
that the value of $\cE_H(\cO)$ should be the same for ``almost all'' orbits $\cO$, in an appropriate sense). 
Except in the case of very simple maps, an exact determination of the Diophantine entropy $\cE_H$ 
is usually at least as difficult as finding the 
algebraic entropy $\cE$, but it  can sometimes be computed explicitly for special families of orbits \cite{markoff}, and a numerical calculation of rational values of iterates along several orbits is much less computationally intensive than computing 
the corresponding sequence of 
rational functions and obtaining their degrees. Empirical evidence suggests that the algebraic entropy and the Diophantine entropy should coincide, so that  $\cE=\cE_H(\varphi)$, but we do not 
know how to prove this except for some particular maps $\varphi$ (cases where $\varphi$ is an integrable map in the Liouville sense  \cite{bruschi, maeda, veselov}, when both entropies are zero, are often 
the most tractable ones). Furthermore, if the dimension $N$ is large then even the exact arithmetic required to compute     $\cE_H(\cO)$ numerically becomes difficult. 

The purpose of this article is to propose that, for discrete dynamical systems that generate Laurent polynomials, and for cluster algebras in particular, the growth of Mahler measure provides another natural notion 
of entropy, which in many cases coincides with the algebraic entropy and the Diophantine entropy, but is much easier to compute numerically. Before proceeding to formulate precise definitions and conjectures, 
we will give some motivation for why restricting to dynamics with the Laurent property need not be unnecessarily restrictive. 

There are many examples of recurrences, birational maps or difference equations 
with the Laurent property. Some of the first known examples, discussed in \cite{gale}, were recurrence relations of the form 
\beq\label{lrec} 
x_{n+N}x_n =F(x_{n+1},\ldots,x_{n+N-1})
\eeq 
for certain polynomials $F$. Notable examples include a recurrence attributed to Dana Scott, that is  
\beq\label{markoff} 
x_{n+3}x_n = x_{n+2}^2+x_{n+1}^2
\eeq 
which produces sequences of Markoff triples:  Markoff's equation 
\beq \label{meq} 
x^2+y^2+z^2=3xyz
\eeq 
arises in Diophantine approximation theory \cite{cassels}, and starting from any triple of positive integers  $(x,y,z)$ satisfying this equation, each subsequent set of three adjacent  
terms generated by (\ref{markoff} is a solution, e.g. taking initial data $(x_1,x_2,x_3)=(1,1,1)$ produces the sequence 
\beq\label{mseq} 
1,1,1,2,5,29,433,37666, 48928105,\ldots 
\eeq 
(A064098 in \cite{oeis}), yielding further Markoff triples $(1,1,2)$, $(1,2,5)$, $(2,5,29)$, etc. Another famous example is the Somos-4 recurrence 
\beq\label{s4}
x_{n+4}x_n = x_{n+3}x_{n+1}+x_{n+2}^2, 
\eeq 
which generates the integer sequence 
\beq\label{s4seq}
1,1,1,1,2,3,7,23,59,314,1529,8209,83313,\ldots  
\eeq 
(A006720  in \cite{oeis}) 
from the initial data $(x_1,x_2,x_3,x_4)=(1,1,1,1)$; for the connection with sequences of points on elliptic curves, see \cite{hones4}. 
In tandem with the development of cluster algebras, an approach to proving the Laurent property for recurrences of the form 
(\ref{lrec}) via the so-called Caterpillar Lemma was presented in \cite{laurentphen}, and it was subsequently shown that a large class of examples with $F$ being a sum of two monomials
can be constructed systematically from  cluster algebra exchange relations (\ref{clustmut}) defined by exchange matrices $B$ with a suitable periodicity property under sequences of 
mutations \cite{FM} (generalized in \cite{nak}), while examples with more than two monomials on the right-hand side arise in the broader setting of LP algebras \cite{lp}. Recurrences with the Laurent property, including 
those of Somos type like (\ref{s4}), 
also occur naturally in the theory of integrable systems, in the form of bilinear discrete Hirota equations for tau functions of lattice equations and their reductions to integrable maps/discrete Painlev\'e 
equations \cite{FH, gal, hi, hkq, in, kun, mase, okubo}.  

However, as pointed out in \cite{conf}, there is a very close connection between the Laurent property and the notion of singularity confinement, which can be combined 
with algebraic entropy as a tool to detect integrability (see \cite{sing} and references), and in this context one finds more general forms of discrete systems with the Laurent property which do not 
seem to arise from iteration of exchange relations $x_k'x_k=F(\bx)$ in a cluster algebra or LP algebra. As described in \cite{hhkq}, given a birational map with confined singularities, which need not have 
the Laurent property, we conjecture 
that it is always possible to construct its ``Laurentification'', that is, a lift to a map in higher dimensions that does have the Laurent property. This conjecture is very natural in the context of discrete integrable systems, 
where one usually  starts from a symplectic map, or a  lattice equation that preserves a Poisson bracket \cite{in}, and then the existence of appropriate tau functions implies that this lifts to a system of discrete Hirota 
equations \cite{zabrodin}, for which  the Laurent property is known to hold \cite{laurentphen}.  Beyond the setting of integrability, Hietarinta and Viallet's example 
\beq\label{hvmap} 
u_{n+1}+u_{n-1}=u_n + \frac{a}{u_n^2}
\eeq 
is a symplectic map of standard type in the plane, with a parameter $a\neq 0$, that has confined singularities but displays chaotic orbits and positive entropy $\cE=\log\big(\frac{1}{2}(3+\sqrt{5})\big)$ \cite{hv,tak}; 
the singularity pattern suggests lifting it to 5 dimensions via 
$
u_n = (x_{n+1}x_{n+2})^{-2}{x_{n+3}x_n}
$,
and this provides a Laurentification of (\ref{hvmap}) because $x_n$ satisfies 
\beq\label{hvl} 
x_{n+5}x_{n+2}^3x_{n+1}^2= 
x_{n+4}^3x_{n+1}^3-x_{n+4}^2 x_{n+3}^3x_{n}+ax_{n+3}^6x_{n+2}^6, 
\eeq 
with $x_n\in\Z[a,x_1^{\pm 1},x_2^{\pm 1},x_3^{\pm 1},x_4^{\pm 1},x_5^{\pm 1}]$ for all $n\in\Z$ \cite{sigmavadim}. For other analogues of tau functions and the Laurent 
property for nonintegrable systems outside the setting of cluster algebras, including multidimensional generalizations of (\ref{hvmap}), see \cite{kkmt} and references.  

We now formulate the main definition of entropy that we would like to consider, for the case of maps $\varphi$ in dimension $N$ with the Laurent property, of the form 
\beq\label{phimap}
\bx \mapsto \varphi (\bx), 
\eeq  
 where  the Laurent property means that for all $n$, each of the components of $\varphi^n(\bx)$ is a Laurent polynomial in the variables appearing in the initial 
data $\bx=(x_1,\ldots ,x_n)$ with integer coefficients. Usually we are interested in birational maps, so this means $\forall n\in\Z$, in which case $\varphi$ can be regarded 
as an automorphism of the field of fractions $\C(x_1,\ldots,x_N)$, but it is also possible to consider non-invertible $\varphi$ and restrict to $n\geq 0$. Then we will 
use $m\big(\varphi^n(\bx)\big)$ to denote the maximum of the Mahler measures of the components of $\varphi^n(\bx)$. 

\begin{definition}\label{EM} For a  map (\ref{phimap}) with the Laurent property, the Mahler entropy $\cE_M$ is 
\beq\label{ME}
\cE_M (\varphi):= \underset {n\to\infty}{\lim} \frac{\log m\big(\varphi^n(\bx)\big)}{n},
\eeq 
whenever the limit exists. 
\end{definition}

\begin{remark}\label{recs} 
All of the examples considered below are of recurrence form, 
$$
\varphi: \, (x_n,\ldots,x_{n+N-1})\mapsto (x_{n+1},\ldots,x_{n+N}), 
$$
with $x_{n+N}$ defined by a relation of the form (\ref{lrec}), so that it is sufficient to take 
\beq\label{alt}
\cE_M:= \underset {n\to\infty}{\lim} \frac{\log m(x_n)}{n},
\eeq 
which turns out to  coincide with the preceding definition in this case. 
\end{remark} 

It is also natural to consider the growth of Mahler measure for a cluster algebra $\cA=\cA(\bx,B)$ defined by an initial seed $(\bx,B)$, by evaluating the Mahler measure $m(x_j')$ for each 
cluster variable $x_j'$ belonging to a cluster $\bx'=\bmu(\bx)$ obtained by a composition of $n$ mutations, and taking the limit $n\to\infty$. In this context, we let $m(\bx')$ denote the maximum 
of the Mahler measures of the elements of the cluster $\bx'$. 

\begin{definition}\label{EMA} For a cluster algebra  $\cA=\cA(\bx,B)$, the Mahler entropy $\cE_M$ is 
\beq\label{MEA}
\cE_M(\cA):= \underset {n\to\infty}{\lim}\underset {|\bmu |=n}{\mathrm{sup}}\, \frac{\log m\big(\bmu(\bx)\big)}{n}, 
\eeq   
assuming that the limit exists. 
\end{definition}

\begin{remark}\label{frozen} 
For cluster algebras with coefficients, given by additional frozen variables that do not mutate, or for maps such as (\ref{hvl}) with parameters, one can either fix specific (integer) values for these coefficients before 
evaluating the Mahler measure, or integrate over an additional torus for each extra parameter that appears, e.g.\ perform an extra integration $\int_{\T^1}  (\cdots) \, \rd \log a$ to evaluate $m(x_n)$ for the case of (\ref{hvl}). 
\end{remark} 

Aside from the issue of existence, it is not obvious that the above definition is independent of the choice of seed, up to mutation equivalence. However, if the Mahler entropy coincides with the 
algebraic entropy, then this is clearly the case. 

\begin{conje}\label{maincon} For cluster algebras $\cA$, or for maps $\varphi$ with the Laurent property obtained from cluster algebras with periodicity (in the sense of \cite{FM} or \cite{nak}), the algebraic and Mahler entropies 
coincide, i.e.\ $\cE=\cE_M$. 
\end{conje} 

\begin{remark}\label{toral}
The assumptions on $\varphi$ mean that it is a composition of finitely many mutations in a cluster algebra, possibly also composed with a  permutation of the coordinates; such maps 
were referred to as cluster maps in \cite{FH}. 
(Several examples will be given below, but see \cite{FH, FM, nak} for many more.)
The above conjecture is certainly false without imposing any restriction on the class of  maps $\varphi$ with the Laurent property, as can be seen by considering the case of monomial maps, which are closely related 
to toral endomorphisms (or automorphisms, if we restrict to birational $\varphi$). Any Laurent monomial $P$ has $m(P)=0$, hence $\cE_M=0$ for all monomial maps, but generic monomial maps have algebraic 
entropy $\cE>0$  \cite{hp}. 
\end{remark}

Henceforth we will mainly be concerned with studying examples that provide evidence for the preceding conjecture. There are various reasons why it seems particularly fruitful to focus on the Mahler entropy for cluster algebras or maps 
generated by compositions of cluster mutations. Firstly, the meromorphic $N$-form (\ref{vol}) is invariant under cluster mutations (\ref{clustmut}) up to a sign, that is $\mu_k^*(\rd \Omega)=\pm \rd \Omega$, which might 
allow for some simplification in explicitly evaluating the Mahler measure of cluster variables. Secondly, it is already interesting to consider Mahler measure in the case of the simplest non-trivial cluster algebra, 
namely the cluster algebra of finite Dynkin type $A_2$,  obtained from the Lyness recurrence 
\beq\label{lyness} 
\varphi: \quad x_{n+2}x_n=x_{n+1}+1, 
\eeq 
which produces the 5-cycle of cluster variables 
\beq\label{5cycle} 
x_1,x_2,x_3=\frac{x_2+1}{x_1}, x_4=\frac{x_1+x_2+1}{x_1x_2},x_5=\frac{x_1+1}{x_2}, 
\eeq 
with the sequence $(x_n)$ repeating with period 5 thereafter. In this case, we  clearly have  $m(x_1)=0=m(x_2)$, and from Jensen's formula it is easy to see that $m(x_3)=0=m(x_5)$, so the only 
non-zero Mahler measure is $m(x_4)=m( x_1+x_2+1)$ which is given by Smyth's formula (\ref{smyth}). (In the latter example it is clear that $\cE_M(\varphi)=\cE_M(\cA)=0=\cE$ because 
the algebra is of finite type.) 
Thirdly, other examples of Mahler measures for cluster variables have explicit formulae in terms of the 
classical dilogarithm $\mathrm{Li}_2$, or the Bloch-Wigner dilogarithm $D$, evaluated at algebraic arguments (see the next section), while the dilogarithm is ubiquitous in the theory of cluster algebras, 
appearing in  the generating functions that produce (part of) mutations as canonical transformations preserving an associated presymplectic structure, and in  terms of the functional identities satisfied by $Li_2$ 
or $D$ \cite{fg,gsv,hi,nak}; 
in particular, the elements of the 5-cycle (\ref{5cycle}), related to the associahedron $K_4$ \cite{fz2}, correspond to the arguments in the five-term relation for $D$ \cite{zagier}. 

Before proceeding to study other concrete examples in the rest of the paper, we present a simple result that connects the Mahler entropy with the Diophantine entropy and the algebraic entropy, but is based 
on the non-trivial fact of positivity for cluster variables: they are Laurent polynomials with positive coefficients, as was proved for the case of skew-symmetric $B$ in \cite{ls}, and for the general case of 
skew-symmetrizable $B$ in \cite{ghkk}. 

\begin{propn}\label{ents} Suppose that $\varphi$ is a cluster map of recurrence type.  Let $\bo = (1,\ldots,1)$ denote the $N$-tuple consisting of $N$ 1s, and let $\cE_H(\bo)$ denote the Diophantine entropy of the orbit
of $\varphi$ with initial data $\bo$. Then (assuming it exists) the Mahler entropy of $\varphi$ satisfies 
\beq\label{lower}
\cE_M\leq \cE_H(\bo).  
\eeq 
Suppose further that $\deg x_n=C\la^n\big(1+o(1)\big)$, $m(x_n)=C'\la^n\big(1+o(1)\big)$ for some real $\la>1$ and $C,C'>0$, so that Conjecture \ref{maincon} holds with 
$\cE =\cE_M=\log \la>0$. 
Then $\cE_H(\bo)=\log \la$ also. 
\end{propn} 

\begin{prf}
For any $P\in\C[\bx]$,  the standard bounds 
$$ 
M(P)\leq L(P) \leq 2^{\mathrm{sdeg}\,P}M(P) 
$$ 
hold, 
where $L$ denotes the length of $P$ (the sum of the absolute values of the coefficients of $P$) and $\mathrm{sdeg}$ denotes the sum of the degrees of $P$ taken in each variable separately. 
Now suppose that $x_n\in\Z_{>0}[x_1^{\pm 1},\ldots,x_N^{\pm 1}]$ is a positive Laurent polynomial generated by a cluster map $\varphi$ of recurrence type. This means that $x_n$ has a canonical 
expression of the form 
\beq\label{lpoly}
x_n = \frac{P_n(\bx)}{\bx^{\bd_n}}  
\eeq 
where the polynomial $P_n\in \Z_{>0}[x_1,\ldots,x_N]$ is not divisible by any $x_j$ for $j\in[1,N]$, and the monomial $\bx^{\bd_n}$ is defined by the integer $N$-tuple $\bd_n$ (d-vector), which consists entirely 
of positive integers for $n$ large enough; see \cite{fziv} or \cite{FH}  for further details. In terms of the total degree of $x_n$ as a rational function of $\bx$ this also means that, for sufficiently large $n$, 
$\deg x_n = \deg P_n$, and clearly $  \mathrm{sdeg}\,P_n\leq N\deg P_n$. Also, by positivity we have $L(P_n)=P_n(\bo)=x_n(\bo)=H\big(x_n(\bo)\big)$, so substituting into the bounds above and taking logs yields 
$$ 
m(x_n)\leq h\big(x_n(\bo)\big)\leq N' \, \deg x_n +m(x_n), \qquad N'=N \log 2.
$$ 
From the upper bound on $m(x_n)$, taking logs and dividing by $n$  immediately gives (\ref{lower}) in the limit $n\to\infty$. On the other hand, if $\deg x_n$ and $m(x_n)$ have the given asymptotic behaviour, 
then substituting this into the same bounds, taking logs and dividing by $n$ gives 
$$ 
\log \la +\frac{\Big(\log C+\log\big(1+o(1)\big)\Big)}{n} \leq 
\frac{\log h\big(x_n(\bo)\big)}{n}  \leq 
\log \la +\frac{\Big(\log C''+\log\big(1+o(1)\big)\Big)}{n}, 
$$ 
where $C''=CN\log 2+C'$, 
and hence  $\cE_H(\bo)=\log \la$ in the limit $n\to\infty$. 
\end{prf}

\begin{remark}\label{diophe} 
With a slightly more detailed analysis, the analogue  of Proposition \ref{ents} should also hold for an arbitrary cluster map given (up to a permutation) by a composition of mutations, and for the corresponding 
entropies for any cluster algebra $\cA$, but a complete proof requires keeping track of the growth of the different components within each cluster. 
\end{remark}

At this stage it is worth pointing out the numerical advantages of working with the Mahler entropy rather than using exact arithmetic to compute growth of degrees or heights of rational numbers. We will describe 
this for a general map $\varphi$ with the Laurent property, but the same considerations apply to computations with cluster variables, with minor modifications. The main idea is to write an approximation 
to the Mahler measure of the  components 
of the $k$th iterate $\varphi^k(\bx) = (x_{1,k},\ldots,x_{N,k})$ 
in terms of a Riemann sum, that is 
\beq\label{sum}
m(x_{j,k})
\approx
\frac{1}{|\cI|}
\sum_{i\in \cI}
\log 
|x_{j,k} (\bx_i^*)|, 
\qquad j \in [1,N],
\eeq 
where the points $\bx_i^*$ for $i$ belonging to some index set $\cI$ are equidistributed on $\T^N$ in the limit $|\cI|\to\infty$. 
There are various equidistribution results in the literature which are relevant to numerical 
approximation of Mahler measure, in particular the work on torsion points in \cite{bir}, which here corresponds to the obvious choice of a regular lattice of points on the flat torus i.e. choosing the components of $\bx_i^*$ to 
be $M$th roots of unity, so that 
$\bx_i^*=(
e^{2\pi\ri j_1/M} ,\ldots,
e^{2\pi\ri j_N/M})$ 
for 
$0\leq j_1,\ldots,j_N\leq M-1$, with 
$|\cI|=M^N$, or the results on Gaussian periods in \cite{habegger}; for a different
approach, based on the doubling map, see \cite{pf}. Another convenient method, which reduces the likelihood of hitting a divergent logarithm when a point on the zero locus is chosen, is to use a Monte Carlo method, 
generating the points $\bx_i^*$ with a (pseudo)random number generator with a uniform distribution on the torus (in some sense, the doubling map used in \cite{pf} can be regarded as a pseudorandom 
number generator).    However,  regardless of what method is used to select equidistributed points $\bx_i^*$, here the main point is that one can regard each choice of point $\bx_i^*$ 
as initial data for $\varphi$, and then one computes its orbit $\bx_i^*,\varphi(\bx_i^*),\ldots,\varphi^n(\bx_i^*)$ up to some 
desired $n$, and for each $k\in [1,n]$ adds the contribution $\log | \cdot |$ of the $j$th component to the sum (\ref{sum}). Thus one obtains an approximation to the Mahler measures of each Laurent polynomial appearing in the 
first $n$ iterates of the orbit, then takes the component with largest Mahler measure in each iterate, and all of the computations are done very rapidly with floating point arithmetic. This is in contrast to the exact arithmetic with rational functions, or with rational numbers, that is required for numerical computation of the algebraic entropy $\cE$, or the 
Diophantine entropy $\cE_H$, respectively. 





\section{Families of Laurent polynomials in rank 2} 
\setcounter{equation}{0}

In this section we consider a family of Laurent polynomials indexed by two positive integers $n,r$, generated by recurrence relations of second order, of the form 
\beq\label{Frec} 
\varphi: \quad
x_{n+2} x_n = F(x_{n+1}), 
\eeq 
where $F$ is specified by 
\beq\label{Fpoly} 
F(x)=x^r+1. 
\eeq 
All recurrences of the form  (\ref{Frec})  preserve the log-canonical symplectic form 
\beq\label{symp} 
\om = \rd \log x_1 \wedge \rd \log x_2, 
\eeq
proportional to (\ref{vol}) when $N=2$, and 
there are more general choices of $F\in\Z[x]$ that give the Laurent property for the map $\varphi$; in particular, it was observed by Speyer that if $F$ is reciprocal, in the sense that 
$F(x)=x^{\deg F}F(x^{-1})$, then (\ref{Frec}) has the Laurent property, but there are other possible choices of $F$ that work (see \cite{conf, numberpoly} for a complete classification).
Here we restrict to the reciprocal    polynomials (\ref{Fpoly}) for $r\geq 1$, because they are precisely what arises from rank 2 cluster algebras $\cA$ defined by the skew-symmetric 
exchange matrices 
$$ 
B=\left(\begin{array}{cc} 
0 & r \\ 
-r & 0 
\end{array} \right), 
$$ 
by composing the mutation 
$\mu_1: \, x_1' x_1 = x_2^r +1$ 
with the permutation $\rho:\, x_1\leftrightarrow x_2$, so that $\varphi=\rho\circ \mu_1=\mu_2\circ \rho$. The periodicity property $\mu_2\circ\mu_1(B)=B$ and the fact that $\mu_2$ is conjugate to $\mu_1$ 
implies that $\varphi$ generates all of the cluster variables, in the form of the sequence $(x_n)_{n\in\Z}$ with $x_n=x_n(x_1,x_2)\in\Z_{>0}[x_1^{\pm 1},x_2^{\pm 1}]$, and the reversibility 
of the recurrence (\ref{Frec}) implies that $x_n(x_1,x_2)=x_{-n+3}(x_2,x_1)$ for all $n$. As a consequence, we have $\cE_M(\cA)=\cE_M(\varphi)$ (evidence of existence will be provided below). 

\begin{figure}
 \centering 
\epsfig{file=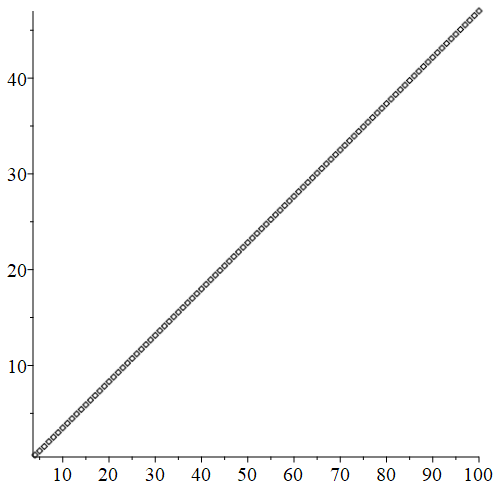, height=2in, width=2in}
\caption{The case  $r=2$: numerical approximation to $m(x_n)$ versus $n$.}\label{sum2} 
\end{figure} 

Similarly, the algebraic entropy $\cE(\cA)$ for the cluster algebra is the same as the algebraic entropy $\cE$ of the map $\varphi$, being given by 
\beq\label{rent} 
\cE=
 \begin{cases}
0  &\text{for}  \,\,r=1,2  \\
 \log\left(\frac{r+\sqrt{r^2-4}}{2}\right)  &\text{for}  \,\, r \geq 3.
\end{cases}
\eeq  
The case $r=1$ is the cluster algebra of finite type $A_2$, generated by the cluster variables (\ref{5cycle}), while for $r=2$ the sequence of polynomials is given explicitly by the formula 
\beq\label{calzel} 
x_{n+3}=\frac{1}{x_1^{n+1}x_2^n}\left(x_2^{2(n+1)}+\sum_{0\leq i+j\leq n}
\Big(\begin{array}{c} n-j \\ i \end{array}\Big)\,\Big(\begin{array}{c} n+1-i \\ j \end{array}\Big)x_1^{2i}x_2^{2j}
\right), 
\eeq 
obtained in \cite{cz}; there is another formula for these polynomials in terms of Chebyshev polynomials \cite{sz}, and below  we will employ yet another closed-form expression for all $n$. 
It is clear from (\ref{calzel}) that the degrees grow linearly with $n$, hence $\cE=0$ for $r=2$. There is an explicit combinatorial formula for the Laurent polynomials $x_n$ in terms of lattice paths \cite{lsank2}, 
valid for any $r$, but for $r\geq 3$ the degree growth is exponential in $n$, and it appears unlikely that there is any  expression for the coefficients as simple as (\ref{calzel}). The entropy value in (\ref{rent}) for 
$r\geq 3$ can be derived directly from the combinatorial formula, or by writing the polynomials in the form (\ref{lpoly}) and using the fact  that the degree vector (d-vector) of the denominator monomial in a cluster algebra satisfies 
the tropical version of the exchange relation, so for the recurrence (\ref{Frec}) with $F$ given by (\ref{Fpoly}) we have 
\beq\label{trop}
\bd_{n+2}+\bd_n = \max (r\bd_{n+1},0), 
\eeq  
satisfied componentwise, with the initial values $\bd_1=(-1,0)^T$, $\bd_2=(0,-1)^T$.  For $n\geq3$, both components of $\bd_n$ are non-negative, so they both satisfy the same scalar 
linear recurrence $d_{n+2}+d_n = rd_{n+1}$, hence $d_n\sim C\la^n$ for some $C>0$, $\la =(r+\sqrt{r^2-4})/2$, and this fixes the same order of growth for the numerator of $x_n$, 
leading to the stated value of the entropy (see e.g.    
\cite{fziv, FH, gal} for more details on tropical dynamics and degree growth).

\begin{figure}
\centering 
\begin{subfigure}{.33\textwidth}
 \centering
 \includegraphics[width=.7\linewidth]{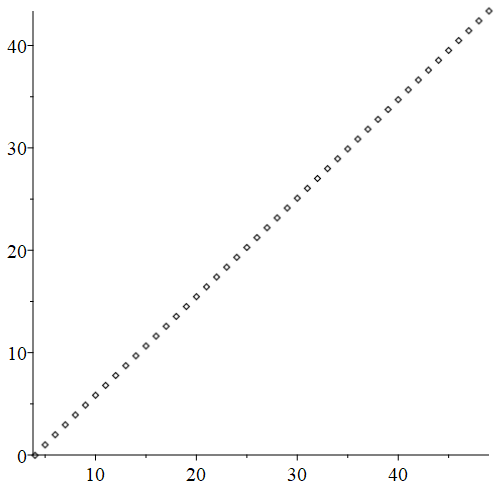}
 \caption{The case  $r=3$.}\label{r3} 
\end{subfigure}%
\begin{subfigure}{.33\textwidth}
 \centering
 \includegraphics[width=.7\linewidth]{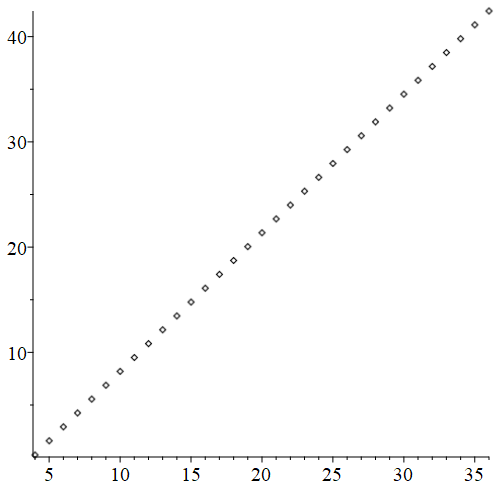}
 \caption{The case  $r=4$.}\label{r4} 
\end{subfigure}%
\begin{subfigure}{.33\textwidth}
 \centering\includegraphics[width=.7\linewidth]{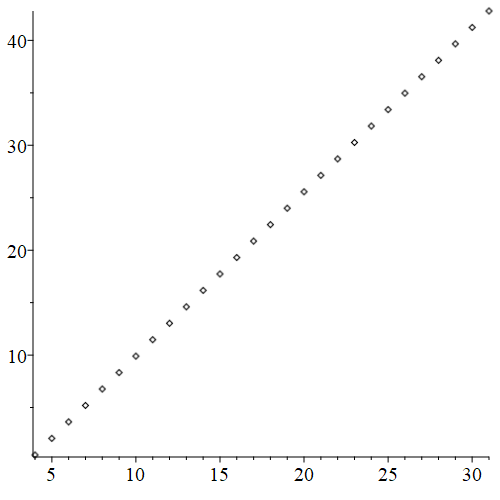}
  \caption{The case  $r=5$.}\label{r5} 
\end{subfigure}\caption{Numerical approximation to $\log m(x_n)$ versus $n$ for $r=3,4,5$.}\label{r45}
\end{figure}

Immediate evidence for Conjecture \ref{maincon} can be seen in Figures \ref{sum2} and \ref{r45}, which were obtained by computing the approximation $S_n\approx m(x_n)$, using (\ref{sum}) with a Monte Carlo method, 
picking $100^2=10^4$ points uniformly distributed on $\T^2$ with a pseudorandom number generator in MAPLE; this took about one minute to produce each plot on a small laptop. For the case $r=2$, the growth 
of $m(x_n)$ appears to be linear in $n$, as is the degree growth of the polynomials  (\ref{calzel}); below we will give a proof of the linear growth of Mahler measure for this case, but for now we just 
record 
the numerical 
estimate 
\beq\label{slope2}
\frac{S_{100}-S_{50}}{50}\approx 0.4837566998
\eeq 
for the slope in Fig.\ref{sum2}, which should be compared with the exact value (\ref{Cstar}) obtained below. 
For $r=3,4,5$ the plots show exponential growth of Mahler measure, and because the growth is so rapid there is a restriction on the size of $n$ that can be calculated before MAPLE gives a floating point error; 
in those cases we have plotted $\log S_n$ against $n$, up to the values $n=49,36,31$, respectively. For $r=3$ we find the estimate 
$\log(S_{49}/S_{25})/{24}\approx                           0.9624236504$ 
for the slope, which should be compared with the value $\cE=\log\big(\tfrac{1}{2}(3+\sqrt{5})\big)\approx 0.9624236498$ for the algebraic entropy in this case. Similarly, for $r=4$ 
we have 
${\log(S_{36}/S_{16})}/{20}\approx    1.316957896$, 
which is extremely close to the algebraic entropy $\cE=\log(2+\sqrt{3})\approx 1.316957897$, and for $r=5$ we find 
${\log(S_{31}/S_{16})}/{15}\approx     1.566799237$, 
which agrees with the value of $\cE=\log\big(\tfrac{1}{2}(5+\sqrt{21})\big)$ to 9 d.p.

We now consider exact evaluation of these Mahler measures. The first two are clearly trivial: $m(x_1)=m(x_2)=0$. More interesting to consider are
the Laurent polynomials 
\beq \label{laur} 
x_3=\frac{x_2^r+1}{x_1}, \quad
x_4 = \frac{(x_2^r+1)^r+x_1^r}{x_1^r x_2}, \quad
x_5 = \frac{1}{x_2^r}\left(\sum_{j=1}^r\Big(\begin{array}{c} r \\ j \end{array}\Big)x_3^{jr-1}+x_1\right), 
\eeq
where (for any $r$) the latter expression for $x_5$ is written most compactly with a sum over powers of $x_3$, making it clear that it is an element 
of $\Z[x_1^{\pm 1},x_2^{\pm 1}]$. Then it is easy to see from Jensen's formula that $m(x_3)=0$. The next two cases will 
require the integral evaluation
\beq\label{lint} 
\int_0^\theta \log |2\sin t| \, \rd t = -\tfrac{1}{2}D(e^{2\ri\theta})=-\frac{1}{2}\sum_{n=1}^\infty \frac{\sin(2n\theta)}{n^2}, \qquad \theta\in\R, 
\eeq
which is used extensively in \cite{lalin}, where $D(z)=\mathrm{Im}\big(\mathrm{Li}_2(z)\big) +\log |z|\arg (1-z)$ in terms 
of the classical dilogarithm $\mathrm{Li}_2(z)=\sum_{n=1}^\infty \tfrac{z^n}{n^2}$ for $|z|<1$ \cite{zagier}. We will also  
need the Chebyshev polynomials of the first kind, defined by $T_r(\cos\theta)=\cos(r\theta)$. 

\begin{thm}
For each positive integer $r$, the Mahler measure of the Laurent polynomial $x_4$ given in (\ref{laur}) is 
\beq\label{mx4} 
m(x_4) = \frac{r}{\pi}D(e^{\frac{\ri\pi}{3}}), 
\eeq 
which equals $r$ times the value (\ref{smyth}) found by Smyth,  while for $r\geq 2$ the Mahler measure of $x_5$ is 
\beq\label{mx5} 
m(x_5) =\frac{r}{\pi}\sum_{j=1}^r (-1)^{j+1}\big(D(e^{\ri rt_j^*}) - D(e^{\ri t_j^*})\big), 
\eeq 
where $0<t_1^*<\cdots <t_r^*<\pi$ are the solutions of 
$2^{r-1}|\sin(\tfrac{rt}{2})|^r=\sin(\tfrac{t}{2})$ on $(0,\pi)$, or equivalently are given by $t_j^*=\arccos X_j^*$, with $X_j^*$
being a root of the polynomial 
\beq\label{Pr} 
P_r(X) = \frac{\big(1-T_r(X)\big)^r}{2^{1-r}(1-X)}-1 
\eeq
on the interval $-1<X<1$.
\end{thm}

\begin{prf} First  note that taking Mahler measures of both sides of (\ref{Frec}) with $F$ given by (\ref{Fpoly}) 
yields the recursion 
\beq\label{mmrec}
m(x_{n+2})+m(x_n) = m(x_{n+1}^r+1),  
\eeq 
so since $m(x_2)=0$ this gives $m(x_4)=m(x_3^r+1) = m\big((x_2^r+1)^r+x_1^r\big)$, as powers of $x_1,x_2$ in the denominator 
make no contribution; this can also be seen directly from the formula in (\ref{laur}). Thus, by applying Jensen's formula to carry out the integration over $x_1$, we 
find 
$$ 
\frac{r}{2\pi\ri}\underset{|x_2^r+1|>1}{\int_{\T^1}} \log |x_2^r+1|\,\rd \log x_2 
= \frac{r}{2\pi\ri}\underset{|y+1|>1}{\int_{\T^1}} \log |y+1|\,\rd \log y=\frac{r}{2\pi\ri}\underset{|2\cos(t/2)|>1}{\int_{-\pi}^\pi} \log |2\cos(t/2)|\,\rd t, 
$$ 
where we substituted $y=x_1^r\in\T^1$ and included an extra factor of $r$ coming from the fact that $y$ winds around the torus $r$ times. 
Then rewriting the integrand as $\log|2\sin s|$ for $s=(t+\pi)/2$, this becomes 
$$
m(x_4)=\frac{r}{\pi}\int_{\frac{\pi}{6}}^{\frac{5\pi}{6}}\log|2\sin s|\,\rd s = -\frac{r}{2\pi}\big(D(e^{\frac{5\ri\pi}{3}})-D(e^{\frac{\ri\pi}{3}})\big), 
$$ by (\ref{lint}), which equals the required answer  (\ref{mx4}) from the fact that  $D(e^{2\ri\theta})$ is an odd function of $\theta$, and for $r=1$ 
this coincides with the value (\ref{smyth}), approximately $0.3230659473$. Similarly, for $x_5$ it is convenient to use (\ref{mmrec}) and $m(x_3)=0$ to 
obtain $$m(x_5)=m(x_4^r+1)=m\big(x_2^{-r}(x_3^r+1)^r+1\big)=m\big((x_3^r+1)^r+x_2^{r}\big)=m\Big(\big(x_1^{-r}(x_2^r+1)^r+1\big)^r+x_2^{r}\Big).$$ 
In the latter expression, we can replace $x_1\to x_1^{-1}$ without changing the Mahler measure, and then substitute $y=x_1^r$ as before (noting the factor 
of $r$ that drops out from the winding number), to find 
$$
m(x_5)=m\Big(\big(y(x_2^r+1)^r+1\big)^r+x_2^r\Big).
$$
As a polynomial in $y$, the argument of $m$ on the right-hand side above is of degree $r$, with roots $y=(x_2^r+1)^{-r}(-1+e^{\ri \pi (2k+1)/r}x_2)$ for $k=0,\ldots, r-1$, so Jensen's 
formula produces a sum of $r$ integrals over suitable values of $x_2\in\T^1$, but each of these integrals transforms to the case $k=0$ by multiplying $x_2$ by a suitable power of 
$e^{2\pi\ri/r}$, giving $r$ copies of the same integral, and changing variables  to $z=e^{\ri\pi/r}x_2$ in this integral yields 
$$ 
m(x_5) = \frac{r}{2\pi\ri}\underset{|\rho(z)|>1}{\int_{\T^1}} \log|\rho(z)|\, \rd\log z, \qquad \rho(z) = \frac{1-z}{(1-z^r)^r}.
$$ 
We have already seen that this vanishes when $r=1$. 
The value of this integral for $r\geq 2$ is given by evaluating the Bloch-Wigner dilogarithm at certain $D$ corresponding to the boundary points where $|\rho(z)|^2=1$ on the torus, given by 
algebraic values of modulus one that are roots of the equation $\rho(z)\bar{\rho}(z^{-1})=1$, which can be converted into a polynomial in $z$.  However, for computational 
purposes it is more effective to set $z=e^{it}$, $t\in(-\pi,\pi]$ and note that $|\rho|=1$ is equivalent to $2^{r-1}|\sin(\tfrac{rt}{2})|^r=|\sin(\tfrac{t}{2})|$ with $t\neq 0$, and by symmetry 
it is sufficient to 
consider the interval $(0,\pi)$, where the modulus sign on the right-hand side can be removed. There are clearly $r$ solutions which we order as $0<t_1^*<\ldots<t_r^*<\pi$, and 
$|\rho|>1$ on the subintervals $[0,t_1^*)$, $(t_2^*,t_3^*)$,..., up to $(t_r^*,\pi]$ for $r$ even, or $(t_{r-1}^*,t_r^*)$ for $r$ odd. For finding the endpoints it is convenient to 
rewrite the equation $|\rho|^{-2}-1$ as $P_r(X)=0$, expressed in the form (\ref{Pr}) with Chebyshev polynomials (note that this is a polynomial of degree    $r^2-1$ because $T_r(1)=0$), 
and then find $t_j^*=\arccos X_j^*$ for $P_r(X_j^*)=0$, $-1<X_j^*<1$. Thus, setting $t_0^*=0$, and $t_{r+1}^*=\pi$ for $r$ even, we arrive at 
$$ 
m(x_5)=\frac{r}{\pi}\underset{|\rho(e^{it})|>1}{\int_0^\pi}\log |\rho(e^{it})|\,\rd t
= \frac{r}{\pi}\sum_{k=0}^{\lc \tfrac{r}{2}\rc}\int_{t_{2k}^*}^{t_{2k+1}^*} \log\frac{|2\sin(\tfrac{t}{2})|}{|2\sin(\tfrac{rt}{2})|^r}\,\rd t, 
$$
and then writing the integrand as $ \log|2\sin(\tfrac{t}{2})|-r\log|2\sin(\tfrac{rt}{2})|$ and applying (\ref{lint}) to each of the latter terms separately, the result (\ref{mx5}) follows. 
\end{prf}

In the rest of this section, we focus on the case $r=2$, which is very special compared with $r\geq 3$: the dynamics of the map $\varphi$ defined by 
$x_{n+2}x_n=x_{n+1}^2+1$ is integrable, since there is a  conserved quantity 
\beq\label{consK} 
K=\frac{x_1}{x_2}+\frac{x_2}{x_1}+\frac{1}{x_1x_2}, \qquad \varphi^*(K)=K, 
\eeq  
whose level sets are conics, and also linearizable, in the sense that the iterates satisfy the linear relation 
\beq\label{linear} 
x_{n+2}-Kx_{n+1}+x_n=0 
\eeq 
on each level set. (For many more examples of cluster variables satisfying linear relations, see \cite{FH,FM,gal,sigmavadim,pyl} and references.)
This means that the general solution can be parametrized explicitly as 
\beq\label{param} 
x_n = \frac{\cos (u+nv)}{\sin v}, 
\eeq (for $K\neq \pm 2$) where the parameters $u,v$ are related to the initial values $x_1,x_2$ by the same formula for $n=1,2$ and the conserved quantity 
is given in terms of $v$ alone by 
$K=2\cos v$.
A short calculation shows that the symplectic form is rewritten in terms of $u,v$ as 
$\rd \log x_1\wedge \rd \log x_2=\rd v\wedge \rd u$, 
so in the case of real dynamics, when $|K|<2$  and there are compact Liouville tori these are the action-angle coordinates, while in the non-compact case $|K|> 2$ one can 
introduce suitable factors of $\ri$ and rewrite (\ref{param}) in terms of hyperbolic functions. (For $K=\pm 2$ the solution of (\ref{linear}) depends linearly on $n$.)   

\begin{thm} For the map $\varphi$ defined by (\ref{Frec}) with (\ref{Fpoly}) for $r=2$, the Mahler entropy is zero, hence also 
$\cE_M(\cA)=0$ for the cluster algebra. More precisely, the Mahler measure of the Laurent polynomials grows linearly with $n$ in this case, that is 
\beq\label{Cstar}
m(x_n) = C^*  n+O(1),  
\quad
with \quad 
C^*= 
\int_{\T^2}\left| \log \Big| \frac{K+\sqrt{K^2-4}}{2}\Big| \right|\,\rd\Omega\approx 0.483997...,  
\eeq 
where $K=K(x_1,x_2)$ is defined by (\ref{consK}). 
\end{thm}

\begin{prf}
Substituting the formula (\ref{param}) into the integral for the Mahler measure gives 
$$
m(x_n) = \int_{\T^2} \log |e^{u+\ri nv}+e^{-(u+\ri nv}| \,\rd \Omega+O(1)=n  \int_{\T^2} |\mathrm{Im}\,v|\,\rd\Omega +O(1),
$$
where $e^{\pm \ri v}$ are the roots of the characteristic quadratic for the linear recurrence (\ref{linear}), and so each is given by 
$\tfrac{1}{2}(K\pm \sqrt{K^2-4})$ for one of the choices of square roots. Hence we have 
$\mathrm{Im}\,v=\pm \log|\tfrac{1}{2}(K+\sqrt{K^2-4})|$, and the result follows. It appears that the numerical evaluation of the integral 
in (\ref{Cstar}) converges much faster with a regular grid of $M^2$ points on $\T^2$ (corresponding to taking $x_1,x_2$ to be $M$th roots of unity) 
than with the Monte Carlo method, and we have checked the value of $C^*$ up to 6 d.p. by taking successively $M=25,50,100,200,500,1000,2000,4000$. 
\end{prf}

\section{Recurrence for Markoff numbers} 
\setcounter{equation}{0}

Given the exchange matrix 
\beq\label{mB}
B=\left(\begin{array}{ccc} 0 & 2 & -2 \\ 
-2 & 0 & 2 \\ 
2 & -2 & 0 \end{array}\right), 
\eeq 
the mutation 
$
\mu_1: \, x_1'x_1=x_2^2+x_3^2
$ 
defines an involution on the cubic surface defined by fixing 
$$ 
K=\frac{x_1^2+x_2^2+x_3^2}{x_1x_2x_3}, 
$$ which is associated with the moduli space of once-punctured 2-tori \cite{cohn}, and corresponds to Markoff's equation (\ref{meq}) when $K=3$. The same is 
true for the mutations $\mu_2,\mu_3$, and by including the cyclic permutation $\rho: \, (x_1,x_2,x_3)\mapsto (x_3,x_1,x_2)$, we find that $\varphi=\rho^{-1}\circ\mu_1$, 
where $\varphi$  corresponds to a single iteration of  (\ref{markoff}), and  $\varphi^2=\rho^{-2}\circ\mu_2\circ\mu_1$ corresponds to two iterations. The exchange matrix 
(\ref{mB}) is cluster mutation-periodic with period 2, in the terminology of \cite{FM}, since  $\mu_2\circ\mu_1(B)=B$. In general, period 2 $B$ matrices lead to a coupling between two different relations, depending on the 
parity of $n$, but 
in this case $\mu_1(B)=-B$, so the exponents on the right-hand side are the same in both mutations, and there is just a  single recurrence. 

It is shown in \cite{hkl} that $\cE=\log\big(\tfrac{1}{2}(1+\sqrt{5})\big)\approx 0.4812118246$ for (\ref{markoff}), and an analogous calculation with the sequence (\ref{mseq}) shows that 
$\cE_H(\bo)$ takes the same value \cite{markoff}.  As is well known \cite{cassels}, certain combinations of involutions of the Markoff surface generate solutions to Pell equations, corresponding to linear dynamics 
and subexponential degree growth, and it turns out that (up to symmetry) the composition of a mutation together with a permutation in $\varphi$ corresponds to the maximum possible growth rate that can be produced by 
mutations alone, so $\cE=\cE(\cA)$ also. To investigate the growth of Mahler measure, we introduce the monomials 
$y_1=\bx^{\bv_1}=x_2/x_1$,  $y_2=\bx^{\bv_2}=x_3/x_2$ 
associated with the integer basis $\bv_1=(-1,1,0)^T$, $\bv_2=(0,-1,1)^T$ for $\mathrm{im}\,B$, leading to the two-dimensional map $\hat{\varphi}$ defined by 
\beq\label{2dmap} 
\hat{\varphi}: \quad y_{n+2}=y_n\left(y_{n+1}+\frac{1}{y_{n+1}}\right), \qquad y_n = x_{n+1}/x_n, 
\eeq 
which is anti-symplectic with respect to the log-canonical 2-form 
\beq\label{ysymp} 
\rd \log y_1\wedge \rd \log y_2 = -\hat{\varphi}^*(\rd \log y_1\wedge \rd \log y_2 ). 
\eeq 
Then for $N=3$ we can make a change of variables to rewrite (\ref{vol}) as
$\rd \Omega = (2\pi\ri)^{-3}\,  \rd \log y_1\wedge \rd \log y_2 \wedge \rd \log x_2$,  and from (\ref{2dmap}) we can recursively calculate 
\beq\label{marec} 
m(x_{n+1}) = m(x_n)+m(y_n)
\eeq 
where the last term above reduces to an integral over the 2-torus, namely 
\beq\label{ymrec} 
m(y_n)=\frac{1}{(2\pi\ri)^2}\int_{\T^2} \log |y_n(y_1,y_2)| \, \rd \log y_1\wedge \rd \log y_2. 
\eeq 
Note that the sequence of $y_n$ are not all Laurent polynomials in $y_1,y_2$, but rational functions in general, because (\ref{2dmap}) does 
not have the Laurent property, so in general $m(y_n)$ is a difference of Mahler measures of polynomials in these two variables.  

To see how this works, note that we have the trivial integrals 
$ 
m(x_1)=m(x_2)=m(x_3)=0
$
and
clearly 
$m(y_1)=m(y_2)=0$. 
Then for 
$x_4=(x_2^2+x_3^2)/{x_{1}}$, 
it follows from (\ref{marec}), (\ref{2dmap}) and (\ref{ymrec})  that 
$m(x_4) = m(y_3) = m\big(y_1(y_2+y_2^{-1})\big) =m(y_2^2+1)=0$ 
by Jensen's formula. Thus the first non-zero value turns out to be the case of 
$x_5 =(x_3^4+x_1^2x_3^2+2x_2^2x_3^2+x_2^4)/({x_2x_1^2})$, 
with 
$$ 
m(x_5)=m(y_4)=m\big(y_2(y_3+y_3^{-1})\big)
=m(y_3^2+1)-m(y_3)
=m\big(y_1^2(y_2+y_2^{-1})^2+1\big).
$$
Then, by  Jensen again, this gives an  integral evaluation equal to twice 
(\ref{smyth}), namely  
\begin{align*}
m(x_5) &= 2\underset{|y_2+y_2^{-1}|<1}{\int_{\T^1}} \log |y_2+y_2^{-1}|^{-1}\, \frac{\rd\log y_2}{2\pi\ri} 
=\frac{2}{\pi}D(e^{\frac{\ri\pi}{3}}).
\end{align*}

\begin{figure}
 \centering 
\epsfig{file=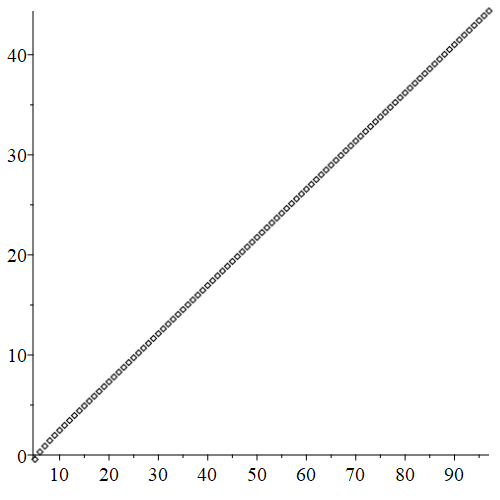, height=2in, width=2in}
\caption{Markoff recurrence: numerical approximation to $\log m(x_n)$ versus $n$ .}\label{summa} 
\end{figure} 

For subsequent terms in the  sequence of Mahler measures $m(y_n)$ given by (\ref{ymrec}), we have done a  numerical evaluation using (\ref{sum}) with a Monte Carlo method, selecting $10^4$ pseudorandom points $(y_1^*,y_2^*)\in \T^2$ as before, iterating the map 
(\ref{2dmap}) for each of these choices of initial data to calculate an approximation to $m(y_n)$,  and 
then applied (\ref{marec}), which allowed us to calculate as far as $S_{97}\approx m(x_{97})$ before getting a floating point error in MAPLE; the results are plotted in Fig.\ref{summa}. The approximate slope is 
$
({\log(S_{97}/S_{50})})/{47}\approx 0.4812118249 
$,
differing from the value of $\cE$ only in the 9th decimal place.

\section{Somos-4 recurrence} 
\setcounter{equation}{0}

Starting from  the exchange matrix 
\beq\label{mBS}
B=\left(\begin{array}{cccc} 0 & 1 & -2 & 1 \\ 
-1 & 0 & 3 & -2 \\ 
2 & -3 & 0 & 1 \\
-1 & 2 & -1 & 0 \end{array}\right), 
\eeq 
the Somos-4 recurrence (\ref{s4}) is given by the composition $\varphi=\rho^{-1}\circ\mu_4\circ\mu_3\circ\mu_2\circ\mu_1$, and this is an example of 
cluster mutation-periodicity with period 1 \cite{FM}. The d-vectors in this case grow quadratically with $n$ \cite{FH}, as do the heights of the integers in (\ref{s4seq}) (this  
is related to the fact that each iteration of the map corresponds to a addition of the same point on an elliptic curve \cite{hones4}), so this is an example with $\cE=0=\cE_H(\bo)$. 
Given a seed defined by (\ref{mBS}), other sequences of mutations produce cluster variables with exponential degree growth, so in this case the algebraic entropy of 
the map $\varphi$ does not equal $\cE(\cA)>0$. 

For an efficient computation of the Mahler measures, we introduce the  coordinates 
\beq\label{ydef} 
y_n=\frac{x_{n+2}x_n}{x_{n+1}^2},  
\eeq 
which satisfy a well-known map of QRT type \cite{qrt}, that is 
\beq\label{qrtmap} 
y_{n+2}=\frac{y_{n+1}+1}{y_ny_{n+1}^2}, 
\eeq 
with the invariant symplectic form $\rd\log y_1\wedge \rd\log y_2$, and again we can rewrite (\ref{vol}), in this case with $N=4$, as a log-canonical form, 
namely 
$ \rd \Omega =  (2\pi\ri)^{-4}\,  \rd \log y_1\wedge \rd \log y_2 \wedge \rd \log x_3 \wedge \rd \log x_4$, 
which 
restricts to the 4-torus  in terms of the new coordinates $(y_1,y_2,x_3,x_4)\in\T^4$. 
Thus from (\ref{ydef}) we obtain the recursion 
\beq\label{smrec} 
m(x_{n+2})=2m(x_{n+1})-m(x_n)+m(y_n). 
\eeq 
Then just as in the Markoff case, we have a sequence of rational functions $y_n(y_1,y_2)$, and the problem of evaluating the Mahler measures 
of the 4-variable Laurent polynomials $x_n(\bx)$ with $\bx=(x_1,x_2,x_3,x_4)$ is reduced to evaluating integrals of the form (\ref{ymrec}) over $\T^2$ and applying a recursion, in this case the relation (\ref{smrec}).   

Now the trivial initial values are 
$m(x_1)=m(x_2)=m(x_3)=m(x_4)=0$ 
and 
$m(y_1)=m(y_2)=0$,
while for 
$y_3=y_1^{-1}y_2^{-2}({y_2+1})$
Jensen's formula gives 
$m(y_3)=0\implies m(x_5)=0$ 
by (\ref{smrec}), where 
$x_5 =x_1^{-1} ({x_4x_2+x_3^2})$.
Thus 
$x_6 =x_1^{-1}x_2^{-1}({x_1x_4^2+x_2x_3x_4+x_3^2})$ 
 is  the first interesting case, 
for which (\ref{smrec}) implies that the Mahler measure is 
$$ 
m(x_6) = m(y_4)=m\left(\frac{y_1y_2(y_1y_2^2+y_2+1)}{(y_2+1)^2}\right)
=m(y_1y_2^2+y_2+1)=\underset{|y_2+1|>1}{\int_{\T^1}}\log |y_2+1|\, \frac{\rd\log y_2}{2\pi\ri},
$$
hence $m(x_6)=\pi^{-1}D(e^{\frac{\ri\pi}{3}})$, the same as the value in Smyth's result (\ref{smyth}). 

\begin{figure}
\centering 
\begin{subfigure}{.33\textwidth}
 \centering
 \includegraphics[width=.7\linewidth]{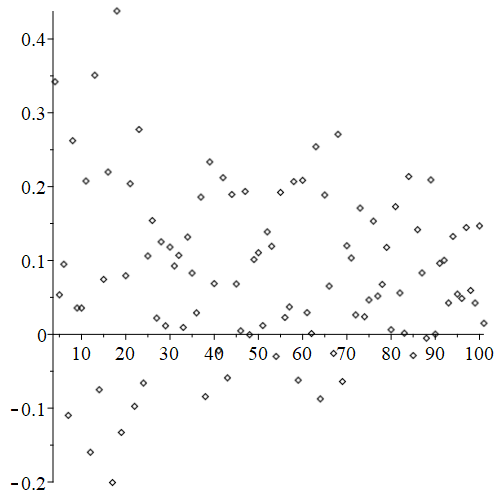}
 \caption{$m(y_n)$ vs. $n$.}\label{s4y} 
\end{subfigure}%
\begin{subfigure}{.33\textwidth}
 \centering
 \includegraphics[width=.7\linewidth]{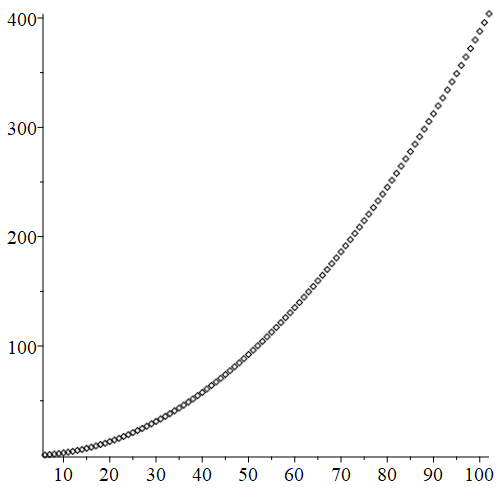}
 \caption{$m(x_n)$ vs. $n$.}\label{s4x} 
\end{subfigure}%
\begin{subfigure}{.33\textwidth}
 \centering\includegraphics[width=.7\linewidth]{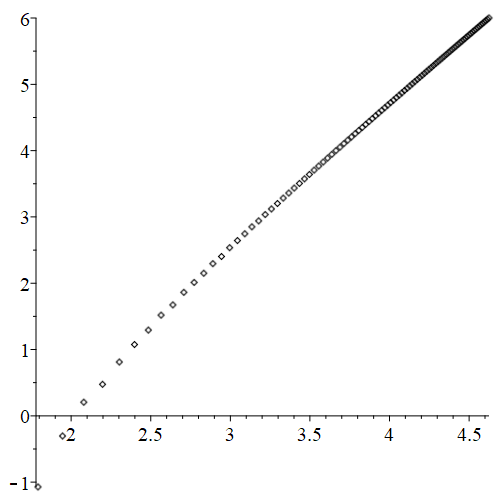}
  \caption{$\log m(x_n)$ vs. $\log n$.}\label{s4logx} 
\end{subfigure}\caption{Numerical approximation for the Somos-4 case.}\label{s4plot}
\end{figure}

Upon doing a numerical calculation of the sequence $m(y_n)$, defined  as integrals over $\T^2$, we observe that  these values seem to oscillate quite erratically, but the corresponding sequence 
of approximations to the Mahler measures $m(x_n)$, obtained from the recursion (\ref{smrec}), does indeed seem to trace out a parabola when plotted against $n$, suggesting that 
\beq\label{quadms4} 
m(x_n) =  Cn^2+O(n), \qquad \mathrm{for}\,\,\,\mathrm{some}\quad C>0, 
\eeq 
which would imply that indeed $\cE_M=0$ in this case. 
This quadratic growth is consistent with the growth of logarithmic heights under iterated translation by a point on an elliptic curve \cite{silverman}, 
and with the growth of degrees for maps of surfaces with an invariant elliptic fibration \cite{df}. 
In Fig.\ref{s4plot} we have also done a log-log plot, which asymptotically looks like a line 
with slope roughly $(\log S_{100}-\log S_{50})/(\log 100 -\log 50)\approx 2.07$, so very close to 2, as expected, and we have estimated the value of $C$ by calculating the 
differences $S_{n+1}-S_n\approx m(x_{n+1})-m(x_n)\sim 2Cn$, which suggests that $C\approx 0.04$.   

To give more support to the conjectured asymptotic formula (\ref{quadms4}), we can refer to the analytic solution to the initial value problem 
for the Somos-4 recurrence, as presented in \cite{hones4}. 
To be precise, the general analytic solution of 
 (\ref{s4}) is given by 
\beq\label{s4sol} 
x_n = A\, B^n \,\frac{\si (z_0+nz)}{\si(z)^{n^2}}, 
\eeq 
where $A,B\in\C^*$ and $\si (\,\cdot\,)=\si (\,\cdot\, ; g_2,g_3)$ is the Weierstrass sigma function for an 
associated cubic curve $\mathrm{y}^2 = 4 \mathrm{x}^3 -g_2\, \mathrm{x} -g_3$. 
In the formula (\ref{s4sol}), the invariants $g_2,g_3\in\C$ of the associated curve and the parameters $z_0,z\in\C$ appearing in the argument of $\si$ depend only on the pair of initial values $y_1,y_2$ for the map 
(\ref{qrtmap}). This implies that the Mahler measure  in this case is given to leading order by an integral over the 2-torus, namely 
\beq\label{s4asy} 
m(x_n) = 
\frac{1}{(2\pi\ri)^2}\int_{\T^2} \Big(\log |\si (z_0 +n z)|-n^2 \,\log|\si(z)| \Big)  \, \rd \log y_1\wedge \rd \log y_2+O(n),  
\eeq 
where $z_0=z_0 (y_1,y_2)$ and $z=z (y_1,y_2)$ are determined by elliptic integrals, and the invariants $g_2,g_3$ are rational functions of $y_1,y_2$. However, as it stands, 
the expression (\ref{s4asy}) does not immediately yield the result  (\ref{quadms4}), because the term 
$\log |\si (z_0 +n z)|$ contains some hidden $n^2$ growth, which is best seen by rewriting this expression in terms of a Jacobi theta function (cf.\ the asymptotic calculation for 
numerical growth of Somos-5 sequences in \cite{hones5}). 

\section{Concluding remarks} 
\setcounter{equation}{0}

It appears that the Mahler entropy gives a useful numerical tool for measuring growth in cluster algebras and dynamical systems with the Laurent property, but so far in rank 2 we have only been able to determine it precisely 
for the simplest examples defined by (\ref{Fpoly}) with $r=1,2$. The numerical results for $r=3,4,5$ suggest that 
$$ 
m(x_{n+2})+m(x_n) =rm(x_{n+1})+O(1), 
$$ 
but it is not clear that this follows immediately from the exact relation 
(\ref{mmrec}). If it should turn out that, to leading order, the sequence of Mahler measures satisfies a tropical version of the original dynamical system, then this could provide a key to 
relating it to the dynamics of d-vectors and degree growth, possibly leading to a proof of Conjecture \ref{maincon}. 

A next step in a more interesting case would be to use the analytic formula  in \cite{hones4}, expressing the solution of  Somos-4 in terms of the Weierstrass 
sigma function, to give a proof of the conjectured asymptotics (\ref{quadms4}), by  deriving an exact integral formula for the constant $C$. So far we have merely indicated 
a possible way to go about this, using the expression (\ref{s4asy}), and we propose to leave the detailed analysis for future work.

It is interesting that all of the Mahler measures of cluster variables that we have been able to evaluate explicitly so far are written in terms of the Bloch-Wigner dilogarithm with algebraic arguments, 
similar to the examples in \cite{lalin}, which are related to volumes of hyperbolic polytopes. Given the connections between cluster algebras and Teichm\"uller theory, it is tempting to 
suggest that there should be a deeper explanation for this phenomenon.

\noindent \textbf{Acknowledgments:} This research was supported by 
Fellowship EP/M004333/1  from the Engineering \& Physical Sciences Research Council, UK, 
and grant 
 IEC\textbackslash R3\textbackslash 193024 from the Royal Society. 
On behalf of all authors, the corresponding author states that there is no conflict of interest.

\end{document}